\documentclass[a4paper,9pt]{article}
\usepackage{amsmath, amsfonts, amssymb}
\usepackage[english, russian]{babel}
\usepackage{graphicx}
\numberwithin{equation}{section}
\begin{document}
\sloppy
\begin{center}
\textbf{On fundamental solutions for multidimensional Helmholtz \\ equation with three singular coefficients}\\[5pt]
\textbf{Ergashev T.G.\\}
\medskip
{ Institute of Mathematics, Uzbek
Academy of Sciences,  Tashkent, Uzbekistan. \\

{\verb ergashev.tukhtasin@gmail.com }}\\
\end{center}

\begin{quote} 

The main result of the present paper is the construction of
fundamental solutions for a class of multidimensional elliptic
equations with three singular coefficients, which could be
expressed in terms of a confluent hypergeometric function of four
variables. In addition, the order of the singularity is determined
and the properties of the found fundamental solutions that are
necessary for solving boundary value problems for
degenerate elliptic equations of second order are found.\\
 \textit{\textbf{Key
words:}} multidimensional elliptic equation with three singular
coefficients, fundamental solutions, confluent hypergeometric
functions of four variables.
\end{quote}

\section{Introduction}

It is known that fundamental solutions have an essential role in
studying partial differential equations. Formulation and solving
of many local and non-local boundary value problems are based on
these solutions. Moreover,\,\,\,fundamental solutions appear as
potentials, for instance, as simple-layer and double-layer
potentials in the theory of potentials.

The explicit form of fundamental solutions gives a possibility to
study the considered equation in detail. For example, in the works
of Barros-Neto and Gelfand \cite{BN}, fundamental solutions for
Tricomi operator, relative to an arbitrary point in the plane were
explicitly calculated. Among other results in this direction, we
note a work by Itagaki \cite{I}, where 3D high-order fundamental
solutions for a modified Helmholtz equation were found. The
fundamental solutions can be applied to some boundary value
problems \cite{{G},{K1},{K2},{KN},{SH8},{SH7}}.

In the present work we find fundamental solutions for the equation
\begin{equation} \label{eq11} {\sum\limits_{i = 1}^{p} {u_{x_{i} x_{i}}} }   +
\frac{{2\alpha _{1}}} {{x_{1}}} u_{x_{1}} +\frac{{2\alpha _{2}}}
{{x_{2}}} u_{x_{2}}+\frac{{2\alpha _{3}}} {{x_{3}}} u_{x_{3}} -
\lambda ^{2}u = 0
\end{equation}
in the domain $R^{3+}_p:=\left\{
\left(x_1,...,x_p\right):x_1>0,x_2>0,x_3>0\right\}\,$, where $p$
is a dimension of the Euclidean space; $p \ge 3$; $\alpha_j$ are
real constants and $0<2\alpha_j<1 $ ($j=1,2,3$); $\lambda$ is real
or pure imaginary constant.

Various modifications of the equation (\ref{eq11}) in the two- and
three-dimensional cases were considered in many papers \cite
{{H},{HK},{UK},{U}}. However, relatively few papers have been
devoted to finding the fundamental solutions for multidimensional
equations, we only note the works \cite{{M},{UK}}. In the paper
\cite{MG}, fundamental solutions of equation (1.1) with a single
singular coefficient are found and investigated.

In this article, at first we shall introduce one confluent
hypergeometric function of four variables. Furthermore, by means
of the introduced hypergeometric function we construct fundamental
solutions of the equation (\ref{eq11}) in an explicit form. For
studying the properties of the fundamental solutions, the
introduced confluent hypergeometric function is expanded in
products of hypergeometric functions of Gauss. With the help of
the obtained expansion it is proved that the constructed
fundamental solutions of equation (\ref{eq11}) have a singularity
of order $1 / r^{p - 2}$ at $r \to 0.$

\section{About one confluent hypergeometric function}

In \cite{E} (see also \cite[p.74,(4b)]{SK}) a hypergeometric
function of many variables of the form
$$
H_{n,p} \left( {a,b_{1} ,...,b_{n} ,c_{p + 1} ,...,c_{n} ;d_{1}
,...,d_{p} ;x_{1} ,x_{2} ,...,x_{n}}  \right) $$
\begin{equation} \label{e21}
= {\sum\limits_{m_{1} ,...,m_{n} = 0}^{\infty}  {}} {\frac{{\left(
{a} \right)_{m_{1} + ... + m_{p} - m_{p + 1} - ... - m_{n}} \left(
{b_{1}} \right)_{m_{1}} ...\left( {b_{n}}  \right)_{m_{n}} \left(
{c_{p + 1}} \right)_{m_{p + 1}}  ...\left( {c_{n}}
\right)_{m_{n}}} } {{\left( {d_{1}} \right)_{m_{1}}  ...\left(
{d_{p}} \right)_{m_{p}}  m_{1} !...m_{n} !}}}x_{1}^{m_{1}}
...x_{n}^{m_{n}},
\end{equation}
 $0 \le p \le n,\,$ is considered,
where $(\kappa)_{m}: = \Gamma (\kappa + m) / \Gamma (\kappa)$ is
the Pochhammer symbol, $m$ is a integer number, $a$ is a complex
number, and $\kappa \neq 0,-1,-2,...,$ if  the Pochhammer symbol
$(\kappa)_m$ is on the denominator.

The hypergeometric function (\ref{e21}) in four variables case has
a form

\begin{equation} \label{e22}H_{4,3} \left( {a,b_1,b_2,b_3, b_4, c_4; d_1,d_2,d_3;
x,y,z,t} \right) = {\sum\limits_{m,n,k,l = 0}^{\infty}  }\frac{
(a)_{m+n+k-l} (b_1)_m(b_2)_n(b_3)_k (b_4)_l(c_4)_l }
{(d_1)_m(d_2)_n(d_3)_km!n!k!l!}x^my^nz^kt^l,
\end{equation}
where $|x|+|y|+|z|<1,\,\, |t|<1/\left(1+|x|+|y|+|z|\right)$.

From the hypergeometric function (\ref{e22}) we shall define the
following confluent hypergeometric function
$$ \textrm{H}_{4,3}^{0} \left({a,b_1,b_2,b_3; d_1,d_2,d_3; x,y,z,t}\right)
 = \mathop {\lim }\limits_{\varepsilon  \to 0} H_{4,3} \left( {a,b_1,b_2,b_3,{\frac{{1}}{{\varepsilon
}}},{\frac{{1}}{{\varepsilon}} };d_1,d_2,d_3;x,y,z,
\varepsilon^{2}t} \right).
$$

At the determination of the hypergeometric function $
\textrm{H}_{4,3}^{0} \left({a,b_1,b_2,b_3; d_1,d_2,d_3;
x,y,z,t}\right)$  the equality  $\lim_{\varepsilon  \to 0}
\left(1/\varepsilon\right)_n\cdot\varepsilon^n=1$ ($n$ is a
natural number) has been used and the found confluent
hypergeometric function represented as
\begin{equation} \label{e23}
 \textrm{H}_{4,3}^{0}
\left({a,b_1,b_2,b_3; d_1,d_2,d_3; x,y,z,t}\right)
 = {\sum\limits_{m,n,k,l = 0}^{\infty}  }\frac{ (a)_{m+n+k-l}
(b_1)_m(b_2)_n(b_3)_k  }
{(d_1)_m(d_2)_n(d_3)_km!n!k!l!}x^my^nz^kt^l, |x|+|y|+|z|<1.
\end{equation}

Using the formula of derivation
$$ \frac{\partial^{i+j+k+l}}{\partial x^i\partial y^j\partial z^k\partial t^l}\textrm{H}_{4,3}^0(a,b_1,b_2,b_3;d_1,d_2,d_3;x,y,z,t)
=\frac{(a)_{i+j+k-l}(b_1)_i(b_2)_j(b_3)_k}{(d_1)_i(d_2)_j(d_3)_k}\times$$
$$\times\textrm{H}_{4,3}^0(a+i+j+k-l,b_1+i,b_2+j,b_3+k;d_1+i,d_2+j,d_3+k;x,y,z,t)
$$ it is easy to show that the confluent hypergeometric function
in (\ref{e23}) satisfies the system of hypergeometric equations
\begin{equation}
\label{e24} {\left\{ {{\begin{array}{*{20}c}
 {\begin{array}{l}
 {x(1-x)\omega_{xx}-xy\omega_{xy}-xz\omega_{xz}+xt\omega_{xt}}
 + {\left[{d_{1}-(a+b_{1}+1)x} \right]}\omega _{x}-
b_{1}y \omega _{y}
-b_1 z \omega _{z}  + b _{1}t \omega _{t}  - ab_{1} \omega = 0 \\
 \end{array}} \hfill \\
 {\begin{array}{l}
 {y (1 - y )\omega _{yy}  - xy \omega _{xy}  - yz \omega _{yz}  + yt \omega
_{yt}}- b_{2} x \omega _{x} + {\left[ {d_{2} - (a + b_{2}+1 )y}
\right]}\omega _{y}
-b_{2} z \omega _{z}  + b_{2} t \omega _{t}  - ab_{2} \omega = 0 \\
 \end{array}} \hfill \\
 {\begin{array}{l}
 {z (1 - z )\omega _{zz}  - xz \omega _{xz}  - yz \omega_{yz}  + zt \omega
_{zt}} - b_{3} x \omega _{x}
-b_{3} y \omega _{y} + {\left[ {d_{3} - (a +b_{3}+1 )z} \right]}\omega _{z} + b_{3} t \omega _{t}  - ab_{3} \omega = 0 \\
 \end{array}} \hfill \\
 {t \omega _{tt}  - x \omega_{xt}  - y \omega _{yt} -z \omega_{zt}+ (1 - a )\omega _{t}  + \omega
 =0},
\hfill \\
\end{array}}}  \right.}
\end{equation}
where $$
\omega(x,y,z,t)=\textrm{H}_{4,3}^0(a,b_1,b_2,b_3;d_1,d_2,d_3;x,y,z,t).$$

Having substituted $\omega(x,y,z,t)=x^\tau y^\nu z^\mu t^\delta
\psi(x,y,z,t)$ in the system of hypergeometric equations
(\ref{e24}), it is possible to find 8 linearly independent
solutions of system (\ref{e24}), which is given in the table form \\
\newcommand{\PreserveBackslash}[1]{\let\temp=\\#1\let\\=\temp}
\let\PBS=\PreserveBackslash
\\
\begin{tabular}
{|p{20pt}|p{26pt}|p{26pt}|p{26pt}|p{26pt}|p{26pt}|p{26pt}|p{26pt}|p{26pt}|}
\hline & $ \omega _{1} $ & $ \omega _{2} $ & $ \omega _{3} $ & $
\omega _{4} $ & $ \omega _{5} $ & $ \omega _{6} $ & $ \omega _{7}
$ & $ \omega _{8} $
 \\
\hline
$
\tau
$
&
$
0
$
&
$
1 - d_{1}
$
&
$
0
$
&
$
0
$
&
$
1 - d_{1}
$
&
$
1 - d_{1}
$
&
$
0
$
&
$
1 - d_{1}
$
 \\
\hline
$
\nu
$
&
$
0
$
&
$
0
$
&
$
1 - d_{2}
$
&
$
0
$
&
$
1 - d_{2}
$
&
$
0
$
&
$
1 - d_{2}
$
&
$
1 - d_{2}
$
 \\
\hline
$
\mu
$
&
$
0
$
&
$
0
$
&
$
0
$
&
$
1 - d_{3}
$
&
$
0
$
&
$
1 - d_{3}
$
&
$
1 - d_{3}
$
&
$
1 - d_{3}
$
 \\
\hline
$
\delta
$
&
$
0
$
&
$
0
$
&
$
0
$
&
$
0
$
&
$
0
$
&
$
0
$
&
$
0
$
&
$
0
$
 \\
\hline

\end{tabular}.
\\
\\
or in explicit form as follows
\begin{equation}\label{e25}
\omega_1(x,y,z,t)=\textrm{H}_{4,3}^0(a,b_1,b_2,b_3;d_1,d_2,d_3;x,y,z,t),
\end{equation}
\begin{equation} \label{e26}
\omega_2(x,y,z,t)=x^{1-d_1}\textrm{H}_{4,3}^0(a+1-d_1,b_1+1-d_1,b_2,b_3;2-d_1,d_2,d_3;x,y,z,t),
\end{equation}
\begin{equation}\label{e27}
\omega_3(x,y,z,t)=y^{1-d_2}\textrm{H}_{4,3}^0(a+1-d_2,b_1,b_2+1-d_2,b_3;d_1,2-d_2,d_3;x,y,z,t),
\end{equation}
\begin{equation}\label{e28}
\omega_4(x,y,z,t)=z^{1-d_3}\textrm{H}_{4,3}^0(a+1-d_3,b_1,b_2,b_3+1-d_3;d_1,d_2,2-d_3;x,y,z,t),
\end{equation}
\begin{equation}\label{e29}
 \omega_5(x,y,z,t)=x^{1-d_1}y^{1-d_2}
\textrm{H}_{4,3}^0(a+2-d_1-d_2, b_1+1-d_1,
b_2+1-d_2,b_3;2-d_1,2-d_2,d_3;x,y,z,t),
\end{equation}
\begin{equation}\label{e210}
\omega_6(x,y,z,t)=x^{1-d_1}z^{1-d_3}
\textrm{H}_{4,3}^0(a+2-d_1-d_3, b_1+1-d_1,
b_2,b_3+1-d_3;2-d_1,d_2,2-d_3;x,y,z,t),
\end{equation}
\begin{equation}\label{e211}
\omega_7(x,y,z,t)=y^{1-d_2}z^{1-d_3}
\textrm{H}_{4,3}^0(a+2-d_2-d_3, b_1,b_2+1-d_2,
b_3+1-d_3;d_1,2-d_2,2-d_3;x,y,z,t),
\end{equation}
\begin{equation}\label{e212}
\omega_8(x,y,z,t)\hfill\\=x^{1-d_1}y^{1-d_2}z^{1-d_3}
\textrm{H}_{4,3}^0(a+3-d_1-d_2-d_3,b_1+1-d_1,b_2+1-d_2,b_3+1-d_3;
2-d_1,2-d_2,2-d_3;x,y,z,t).
\end{equation}

\bigskip
\section{Decomposition formulas}

For a given multivariable function, it is useful to find a
decomposition formula which would express the multivariable
function in terms of products of several simpler hypergeometric
functions involving fewer variables. For this purpose Burchnall
and Chaundy \cite{BC} had given a number of expansions of double
hypergeometric functions in series of simpler hypergeometric
functions. Their method is based upon the inverse pair of symbolic
operators
\begin{equation}
\label{eq31} \nabla \left( {h} \right): = {\frac{{\Gamma \left(
{h} \right)\Gamma \left( {\delta _{1} + \delta _{2} + h}
\right)}}{{\Gamma \left( {\delta _{1} + h} \right)\Gamma \left(
{\delta _{2} + h} \right)}}}, \quad \Delta \left( {h} \right): =
{\frac{{\Gamma \left( {\delta _{1} + h} \right)\Gamma \left(
{\delta _{2} + h} \right)}}{{\Gamma \left( {h} \right)\Gamma
\left( {\delta _{1} + \delta _{2} + h} \right)}}},
\end{equation}
where $$\delta _{1} : = x_{1} {\frac{{\partial}} {{\partial
x_{1}}} }, \delta _{2} : = x_{2} {\frac{{\partial}} {{\partial
x_{2}}} }$$

Recently Hasanov and Srivastava \cite {{HS6},{HS7}} generalized
the operators $\nabla \left( {h} \right)$ and $\Delta \left( {h}
\right)$ defined by (\ref{eq31}) in the forms
\begin{equation}
\label{eq32} \tilde {\nabla} _{x_{1} ;x_{2} ,...,x_{m}}  \left(
{h} \right): = {\frac{{\Gamma \left( {h} \right)\Gamma \left(
{\delta _{1} + ... + \delta _{m} + h} \right)}}{{\Gamma \left(
{\delta _{1} + h} \right)\Gamma \left( {\delta _{2} + ... + \delta
_{m} + h} \right)}}}
\end{equation}
and
\begin{equation}
\label{eq33} \tilde {\Delta} _{x_{1} ;x_{2} ,...,x_{m}}  \left(
{h} \right): = {\frac{{\Gamma \left( {\delta _{1} + h}
\right)\Gamma \left( {\delta _{2} + ... + \delta _{m} + h}
\right)}}{{\Gamma \left( {h} \right)\Gamma \left( {\delta _{1} +
... + \delta _{m} + h} \right)}}},
\end{equation}
where \begin{equation}\label{eq331}\delta _{i} : = x_{i}
{\frac{{\partial}} {{\partial x_{i}}} }\,\,\,\left( {i = 1,...,m}
\right),
\end{equation}
and they obtained very interesting results. For example, in the
special case when $m = 3$ a Lauricella function in three variables
is defined by (cf.\cite{L}; see also \cite[p.33, 1.4(1)]{SK})
$$
F_{A}^{\left( {3} \right)} \left( {a,b_{1} ,b_{2} ,b_{3} ;d_{1}
,d_{2} ,d_{3} ;x,y,z} \right) = {\sum\limits_{l,m,n = 0}^{\infty}
{{\frac{{\left( {a} \right)_{l + m + n} \left( {b_{1}}
\right)_{l} \left( {b_{2}} \right)_{m} \left( {b_{1}}
\right)_{n}}} {{\left( {d_{1}}  \right)_{l} \left( {d_{2}}
\right)_{m} \left( {d_{1}}  \right)_{n}}} }}
}{\frac{{x^{l}}}{{l!}}}{\frac{{y^{m}}}{{m!}}}{\frac{{z^{n}}}{{n!}}}
$$
and the following decomposition formula holds true \cite{HS6}
\begin{equation}
\label{eq34}
\begin{array}{l}
 F_{A}^{(3)} \left( {a;b_{1} ,b_{2} ,b_{3} ;d_{1} ,d_{2} ,d_{3} ;x,y,z}
\right) = {\sum\limits_{l,m,n = 0}^{\infty}  {{\frac{{\left( {a}
\right)_{l + m + n} \left( {b_{1}}  \right)_{l + m} \left( {b_{2}}
\right)_{l + n} \left( {b_{3}}  \right)_{m + n}}} {{\left( {d_{1}}
\right)_{l + m} \left( {d_{2}}  \right)_{l + n} \left( {d_{3}}
\right)_{m + n} l!m!n!}}}x^{l +
m}y^{l + n}z^{m + n}}}  \\
 \cdot F\left( {a + l + m,b_{1} + l + m;d_{1} + l + m;x} \right)F\left( {a +
l + m + n,b_{2} + l + n;d_{2} + l + n;y} \right) \\
 \cdot F\left( {a + l + m + n,b_{3} + m + n;d_{3} + m + n;z} \right), \\
 \end{array}
\end{equation}
where $$F\left( {a,b;c;x} \right) = {\sum\limits_{n = 0}^{\infty}
{{\frac{{\left( {a} \right)_{n} \left( {b} \right)_{n}}} {{\left(
{c} \right)_{n} n!}}}x^{n}}} $$ is a Gaussian hypergeometric
function \cite [p.56,(2)]{E}.

It should be noted that the symbolic notations (\ref{eq331}) in
the one-dimensional case take the form $\delta : = xd / dx$ and
such a notation is used in solving problems of the operational
calculus \cite[p.26]{P}.

We now introduce here the other multivariable analogues of the
Burchnall-Chaundy symbolic operators $\nabla \left( {h} \right)$
and $\Delta \left( {h} \right)$ defined by (\ref{eq31}):
\begin{equation}
\label{eq35} \tilde {\nabla} _{x,y}^{m,n} \left( {h} \right): =
{\frac{{\Gamma \left( {h} \right)\Gamma \left( {h + \delta _{1} +
... + \delta _{m} - \sigma _{1} - ... - \sigma _{n}}
\right)}}{{\Gamma \left( {h + \delta _{1} + ... + \delta _{m}}
\right)\Gamma \left( {h - \sigma _{1} - ... - \sigma _{n}}
\right)}}} = {\sum\limits_{s = 0}^{\infty}  {{\frac{{\left( { -
\delta _{1} - ... - \delta _{m}}  \right)_{s} \left( {\sigma _{1}
+ ... + \sigma _{n}} \right)_{s}}} {{(h)_{s} s!}}}}} ,
\end{equation}

\begin{equation}
\label{eq36} \tilde {\Delta} _{x,y}^{m,n} \left( {h} \right): =
{\frac{{\Gamma \left( {h + \delta _{1} + ... + \delta _{m}}
\right)\Gamma \left( {h - \sigma _{1} - ... - \sigma _{n}}
\right)}}{{\Gamma \left( {h} \right)\Gamma \left( {h + \delta _{1}
+ ... + \delta _{m} - \sigma _{1} - ... - \sigma _{n}} \right)}}}
= {\sum\limits_{s = 0}^{\infty}  {{\frac{{\left( {\delta _{1} +
... + \delta _{m}}  \right)_{s} \left( { - \sigma _{1} - ... -
\sigma _{n}} \right)_{s}}} {{(1 - h)_{s} s!}}}}} ,
\end{equation}
where
\begin{equation}
\label{eq37} x: = \left( {x_{1} ,...,x_{m}}  \right), y: = \left(
{y_{1} ,...,y_{n}}  \right),
\end{equation}
$$ \delta _{i} : = x_{i}
{\frac{{\partial}} {{\partial x_{i}}} }\,\,, \sigma _{j} : = y_{j}
{\frac{{\partial}} {{\partial y_{j}}} }, \quad i = 1,...,m,\,\,j =
1,...,n;\,\,m,n \in {\rm N}. $$

In addition, we consider operators which coincide with Hasanov-
Srivastava's symbolic operators $\tilde {\nabla} \left( {h}
\right)$ and $\tilde {\Delta} \left( {h} \right)$ defined by
(\ref{eq32}) and (\ref{eq33}) as a particular case: $$ \tilde
{\nabla} _{x, -} ^{m,0} \left( {h} \right): = \tilde {\nabla}
_{x_{1} ;x_{2} ,...,x_{m}} \left( {h} \right), \quad \tilde
{\Delta} _{x, -} ^{m,0} \left( {h} \right): = \tilde {\Delta}
_{x_{1} ;x_{2} ,...,x_{m}}  \left( {h} \right), m \in {\rm N}; $$

$$ \tilde {\nabla} _{ - ,y}^{0,n}
\left( {h} \right): = \tilde {\nabla} _{ - y_{1} ; - y_{2} ,..., -
y_{n}}  \left( {h} \right), \quad \tilde {\Delta} _{ - ,y}^{0,n}
\left( {h} \right): = \tilde {\Delta} _{ - y_{1} ; - y_{2} ,..., -
y_{n}}  \left( {h} \right), n \in {\rm N}. $$

It is obvious that
$$
\tilde {\nabla} _{x, -} ^{1,0} \left( {h} \right) = \tilde
{\Delta} _{x, - }^{1,0} \left( {h} \right) = \tilde {\nabla} _{ -
,y}^{0,1} \left( {h} \right) = \tilde {\Delta} _{ - ,y}^{0,1}
\left( {h} \right) = 1.
$$
We introduce the notation:
$$D_z^{s}f(z)={\sum\limits_{I(k,s)} {{\frac{{z_{1}^{i_{1}}
...z_{k}^{i_{k}}} } {{i_{1} !...i_{k} !}}}{\frac{{\partial
^{s}f}}{{\partial z_{1}^{i_{1}}  ...\partial z_{k}^{i_{k}}} }
}}},$$ where
$$z=(z_1,...,z_k); \,\,\,\,I(k,s)=\left\{\left(i_1,...,i_k\right):i_{1} \ge 0,...,i_{k}
\ge 0,i_{1} + ... + i_{k} = s \right\}.$$

\textbf{Lemma.} Let be $f: = f\left( {x} \right)$ and $g: =
g\left( {y} \right)$ functions with variables $x$ and $y$ in
(\ref{eq37}). Then following equalities hold true for any $m,\,n
\in {\rm N}$:
\begin{equation}
\label{eq401} \left( { - x_{1} {\frac{{\partial}} {{\partial
x_{1}}} } - ... - x_{m} {\frac{{\partial}} {{\partial x_{m}}} }}
\right)_{s} f(x)= ( - 1)^{s}s!D_x^{s}f(x) ,\,\,\,s \in {\rm N}
\cup {\left\{ {0} \right\}};
\end{equation}
\begin{equation}
\label{eq12} \left( {y_{1} {\frac{{\partial}} {{\partial y_{1}}} }
+ ... + y_{n} {\frac{{\partial}} {{\partial y_{n}}} }} \right)_{s}
g(y) = {\left\{ {{\begin{array}{*{20}c}
{g(y),\,\,\,\,\,\,\,\,\,\,\,\,\,\,\,\,\,\,\,\,\,\,\,\,\,\,\,\,\,\,\,\,\,\,\,\,\,\,\,\,\,\,\,\,\,\,\,\,\,\,\,\,\,\,\,\,\,
if\,\,\,s=0,} \hfill \\
 {s!{\sum\limits_{q = 1}^{s} {C_{s}^{q} C_{s - 1}^{q - 1} \cdot \left( {s -
q} \right)!}} D_y^sg(y),\,\, if \,\,\,s \in {\rm N}.} \hfill \\
\end{array}}}  \right.}
\end{equation}

The lemma is proved by method of mathematical induction \cite{ET}.

In the present paper we shall use two particular cases ($m =
3$\,\,and\,$n = 1)$\, of the formulas (\ref{eq401}) and
(\ref{eq12}):

\begin{equation}
\label{eq13} \left( { - x{\frac{{\partial}} {{\partial x}}} -
y{\frac{{\partial }}{{\partial y}}} - z{\frac{{\partial}}
{{\partial z}}}} \right)_{s} f(x,y,z) = ( -
1)^{s}s!D^s_{x,y,z}f(x,y,z) ,\,\,\,s \in {\rm N} \cup {\left\{ {0}
\right\}};
\end{equation}

\begin{equation}
\label{eq14} \left( {t{\frac{{d}} {{dt}}}} \right)_{s} g(t) =
{\left\{ {{\begin{array}{*{20}c}
{g(t),\,\,\,\,\,\,\,\,\,\,\,\,\,\,\,\,\,\,\,\,\,\,\,\,\,\,\,\,\,\,\,\,\,\,\,\,\,\,\,\,\,\,\,\,\,\,\,\,\,\,\,\,\,\,\,\,\,s
= 0,} \hfill \\
 {{\sum\limits_{q = 1}^{s} {C_{s}^{q} C_{s - 1}^{q - 1} \cdot \left( {s - q}
\right)!t^{q}g^{(q)}(t)}} ,\,\,\,\,s \in
{\rm N}.} \hfill \\
\end{array}}}  \right.}
\end{equation}

Using the formulas (\ref{eq35}) and (\ref{eq36}), we obtain

\begin{equation}
\label{eq15} {\rm H}_{4,3}^{0} \left( {a;b_{1} ,b_{2} ,b_{3}
;d_{1} ,d_{2} ,d_{3} ;x,y,z,t} \right) = \tilde {\nabla}
_{x,y,z,t}^{3,1} \left( {a} \right)F_{A}^{(3)} \left( {a;b_{1}
,b_{2} ,b_{3} ;d_{1} ,d_{2} ,d_{3} ;x,y,z} \right){}_{0}F_{1}
\left( {1 - a; - t} \right),
\end{equation}
where ${}_{0}F_{1} \left( {a;x} \right) = {\sum\limits_{n =
0}^{\infty} {{\frac{{x^{n}}}{{\left( {a} \right)_{n} n!}}}}}$ is a
generalized hypergeometric function \cite[Chapter IV]{EM}.

Now considering the equalities   (\ref{eq35}), (\ref{eq13}) and
(\ref{eq14}) from the formula (\ref{eq15}) we have
\begin{equation}
\label{eq16}
\begin{array}{l}
 {\rm H}_{4,3}^{0} \left( {a;b_{1} ,b_{2} ,b_{3} ;d_{1} ,d_{2} ,d_{3} ;x_{1}
,x_{2} ,x_{3} ,y} \right)
 = {\sum\limits_{s = 0}^{\infty}  {{\sum\limits_{q = 0}^{s}
{{\sum\limits_{I(3,s)} {A\left( {s,q} \right)C_{s}^{q} {\frac{{( -
1)^{s + q}\left( {b_{1}}  \right)_{i} \left( {b_{2}}  \right)_{j}
\left( {b_{3}} \right)_{k} }}{{\left( {1 - a} \right)_{q} \left(
{d_{1}} \right)_{i} \left( {d_{2}} \right)_{j} \left( {d_{3}}
\right)_{k}}} }}}} }}
}{\frac{{x^{i}}}{{i!}}}{\frac{{y^{j}}}{{j!}}}{\frac{{z^{k}}}{{k!}}}{\frac{{t^{q}}}{{q!}}}
\\
 \cdot F_{A}^{(3)} \left( {a + s;b_{1} + i,b_{2} + j,b_{3} + k;d_{1} +
i,d_{2} + j,d_{3} + k;x,y,z} \right){}_{0}F_{1} \left( {1 - a + q;
- t}
\right), \\
 \end{array}
\end{equation}
where $$A\left( {s,q} \right) = {\left\{ {{\begin{array}{*{20}c}
 {1,\,\,\,{\rm i}{\rm f}\,\,s = 0\,\,{\rm a}{\rm n}{\rm d}\,\,q = 0,} \hfill
\\
 {q / s,\,\,{\rm i}{\rm f}\,\,s \ge 1\,\,{\rm a}{\rm n}{\rm d}\,\,q \ge 0,}
\hfill \\
\end{array}}}  \right.}  I(3,s)=\left\{(i,j,k):i \ge 0,j\ge 0,j
\ge 0,i+j+k=s \right\}.$$

Applying the decomposition formula (\ref{eq34}) to the expansion
(\ref{eq16}), we obtain

\begin{equation}
\label{eq17}
\begin{array}{l}
 {\rm H}_{4,3}^{0} \left( {a;b_{1} ,b_{2} ,b_{3} ;d_{1} ,d_{2} ,d_{3}
;x,y,z,t} \right) \\
 = {\sum\limits_{l,m,n,s = 0}^{\infty}  {{\sum\limits_{q = 0}^{s}
{{\sum\limits_{I(3,s)} {A\left( {s,q} \right)s!C_{s}^{q} {\frac{{(
- 1)^{s + q}\left( {a} \right)_{l + m + n + s} \left( {b_{1}}
\right)_{i + l + m} \left( {b_{2}}  \right)_{j + l + n} \left(
{b_{3}}  \right)_{k + m + n}}} {{\left( {a} \right)_{s} \left( {1
- a} \right)_{q} \left( {d_{1}}  \right)_{i + l + m} \left(
{d_{2}}  \right)_{j + l + n} \left( {d_{3}}  \right)_{k + m + n}
}}}}}} }} } {\frac{{x^{i + l + m}}}{{i!n!}}}{\frac{{y^{j + l +
n}}}{{j!m!}}}{\frac{{z^{k + m + n}}}{{k!l!}}}{\frac{{t^{q}}}{{q!}}} \\
 \cdot F\left( {a + s + l + m,b_{1} + i + l + m;d_{1} + i + l + m;x}
\right)F\left( {a + s + l + m + n,b_{2} + j + l + n;d_{2} + j + l
+ n;y}
\right) \\
 \cdot F\left( {a + s + l + m + n,b_{3} + k + m + n;d_{3} + k + m + n;z}
\right){}_{0}F_{1} \left( {1 - a + q; - t} \right). \\
 \end{array}
\end{equation}

By virtue of the formula \cite[p.64,(22)] {EM}
$$
F\left( {a,b;c;x} \right) = \left( {1 - x} \right)^{ - b}F\left(
{c - a,b;c;{\frac{{x}}{{x - 1}}}} \right),
$$
the expansion (\ref{eq17}) yields

\begin{equation} \label{eq18}
\begin{array}{l}
 {\rm H}_{4,3}^{0} \left( {a;b_{1} ,b_{2} ,b_{3} ;d_{1} ,d_{2} ,d_{3}
;x,y,z,t} \right) = \left( {1 - x} \right)^{ - b_{1}} \left( {1 -
y}\right)^{ - b_{2}} \left( {1 - z} \right)^{ - b_{3}}  \\
 \cdot {\sum\limits_{l,m,n,s = 0}^{\infty}  {{\sum\limits_{q = 0}^{s}
{{\sum\limits_{I(3,s)} {A\left( {s,q} \right)s!C_{s}^{q} {\frac{{(
- 1)^{s + q}\left( {a} \right)_{l + m + n + s} \left( {b_{1}}
\right)_{i + l + m} \left( {b_{2}}  \right)_{j + l + n} \left(
{b_{3}}  \right)_{k + m + n}}} {{\left( {a} \right)_{s} \left( {1
- a} \right)_{q} \left( {d_{1}}  \right)_{i + l + m} \left(
{d_{2}}  \right)_{j + l + n} \left( {d_{3}}  \right)_{k + m + n}
i!j!k!n!m!l!q!}}}}}} }} }t^q  \\
 \cdot\left( {{\frac{{x}}{{1 - x}}}} \right)^{i + l + m}F\left( {d_{1} - a - j - k,b_{1} + i
+ l + m;d_{1} +
i + l + m;{\frac{{x}}{{x - 1}}}} \right) \\
 \cdot \left( {{\frac{{y}}{{1 -
y}}}} \right)^{j + l + n}F\left( {d_{2} - a - i - k - m,b_{2} + j
+ l + n;d_{2} + j + l +
n;{\frac{{y}}{{y - 1}}}} \right) \\
 \cdot \left( {{\frac{{z}}{{1 - z}}}} \right)^{k
+ m + n} F\left( {d_{3} - a - i - j - l,b_{3} + k + m + n;d_{3} +
k + m + n;{\frac{{z}}{{z - 1}}}} \right){}_{0}F_{1} \left( {1 - a
+ q; - t} \right).
\\
 \end{array}
\end{equation}

Expansion (\ref{eq18}) will be used for studying properties of the
fundamental solutions.

\section{Fundamental solutions}

We consider equation (\ref{eq11}) in $R_{p}^{3 +}.$ Let
$x:=(x_1,...,x_p)$ be any point and $x_{0}:=(x_{01},...,x_{0p}) $
be any fixed point of $R_{p}^{3 +}.$ We search for a solution of
(\ref{eq11}) as follows: \begin{equation} \label{eq41} u(x,x_0) =
P(r)w(\sigma ),
\end{equation}
where
$$
\sigma = \left( {\sigma _{1} ,\sigma _{2} ,\sigma _{3} ,\sigma
_{4}} \right),\,\,r^{2} = {\sum\limits_{i = 1}^{p} {(x_{i} -
x_{0i} )^{2}}} , \quad r_{k}^{2} = (x_{k} + x_{0k} )^{2} +
{\sum\limits_{i = 1,i \ne k}^{p} {(x_{i} - x_{0i} )^{2}}},
\,\,P(r) = \left( r^{2} \right)^{ - \alpha}, \,\,\,$$ $$ \alpha =
\alpha_1+\alpha_2+\alpha_3 - 1 + {\frac{{p}}{{2}}}, \,\, \sigma
_{k} = {\frac{{r^{2} - r_{k}^{2}}} {{r^{2}}}} = - {\frac{{4x_{k}
x_{0k}}} {{r^{2}}}}, \quad k = 1,2,3; \quad \sigma _{4} =
-\frac{1}{4}\lambda ^2r^2.
$$

We calculate all necessary derivatives and substitute them into
equation (\ref{eq11}):
\begin{equation}\label{eq42}
{\sum\limits_{m = 1}^{4} {A_{m} \omega _{\sigma _{m} \sigma _{m}}}
}   + {\sum\limits_{m = 1}^{3} {{\sum\limits_{n = m + 1}^{3}
{B_{m,n} \omega _{\sigma _{m} \sigma _{n}}} } } }  +
{\sum\limits_{m = 1}^{3} {C_{m} \omega _{\sigma _{m} \sigma _{4}}}
}   + {\sum\limits_{m = 1}^{4} {D_{m} \omega _{\sigma _{m}}} }   +
E\omega = 0,
\end{equation}where
$$ A_{k} = - {\frac{{4P(r)}}{{r^{2}}}}{\frac{{x_{k}}} {{x_{0k}}}
}\sigma _{k} (1 - \sigma _{k} ), C_{k} =
{\frac{{4P(r)}}{{r^{2}}}}{\frac{{x_{0k}}} {{x_{k}}} }\sigma _{k}
\sigma _{4} + {\frac{{\lambda ^{2}}}{{2}}}P(r)\sigma _{k} , $$  $$
B_{k,l} = {\frac{{4P(r)}}{{r^{2}}}}\left( {{\frac{{x_{0k}}}
{{x_{k}}} } + {\frac{{x_{0l}}} {{x_{l}}} }} \right)\sigma _{k}
\sigma _{l} ,\,\,\,k \ne l, l = 1,2,3, $$
$$ D_{k} = - {\frac{{4P(r)}}{{r^{2}}}}{\left\{ { - \sigma _{k}
{\sum\limits_{m = 1}^{3} {{\frac{{x_{0m}}} {{x_{m}}} }\alpha
_{m}}}   + {\frac{{x_{0k} }}{{x_{k}}} }{\left[ {2\alpha _{k} -
\alpha \sigma _{k}}  \right]}} \right\}},\,\,\quad A_{4} = \lambda
^{2}P(r)\sigma _{4} ,\,\,$$
$$D_{4} = {\frac{{4P(r)}}{{r^{2}}}}\sigma
_{4} {\sum\limits_{m = 1}^{3} {{\frac{{x_{0m}}} {{x_{m}}} }\alpha
_{m} + \lambda ^{2}P(r)\alpha}},  E = {\frac{{4\alpha
P(r)}}{{r^{2}}}}{\sum\limits_{m = 1}^{3} {{\frac{{x_{0m}}}
{{x_{m}}} }-\lambda^2P(r)}}.$$

Using the above given representations of coefficients we simplify
equation (\ref{eq42}) and obtain the following system of
equations:
\begin{equation}
 \label{eq43} {\left\{ {{\begin{array}{*{20}c}
 {\begin{array}{l}
 {\sigma _{1} (1 - \sigma _{1} )\omega _{\sigma _{1} \sigma _{1}}  - \sigma
_{1} \sigma _{2} \omega _{\sigma _{1} \sigma _{2}}  - \sigma _{1}
\sigma _{3} \omega _{\sigma _{1} \sigma _{3}}  + \sigma _{1}
\sigma _{4} \omega
_{\sigma _{1} \sigma _{4}}}   \\
 { + {\left[ {2\alpha_{1} - (\alpha+\alpha_1+1)\sigma _{1}}
\right]}\omega _{\sigma _{1}}  - \alpha _{1} \sigma _{2} \omega
_{\sigma _{2}}  - \alpha _{1} \sigma _{3} \omega _{\sigma _{3}}  +
\alpha _{1} \sigma
_{4} \omega _{\sigma _{4}}  - \alpha \alpha _{1} \omega = 0} \\
 \end{array}} \hfill \\
 {\begin{array}{l}
 {\sigma _{2} (1 - \sigma _{2} )\omega _{\sigma _{2} \sigma _{2}}  - \sigma
_{1} \sigma _{2} \omega _{\sigma _{1} \sigma _{2}}  - \sigma _{2}
\sigma _{3} \omega _{\sigma _{2} \sigma _{3}}  + \sigma _{2}
\sigma _{4} \omega
_{\sigma _{2} \sigma _{4}}}   \\
 { + {\left[ {2\alpha _{2} - (\alpha + \alpha _{2} +1)\sigma _{2}}
\right]}\omega _{\sigma _{2}}  - \alpha _{2} \sigma _{1} \omega
_{\sigma _{1}}  - \alpha _{2} \sigma _{3} \omega _{\sigma _{3}}  +
\alpha _{2} \sigma
_{4} \omega _{\sigma _{4}}  - \alpha \alpha _{2} \omega = 0} \\
 \end{array}} \hfill \\
 {\begin{array}{l}
 {\sigma _{3} (1 - \sigma _{3} )\omega _{\sigma _{3} \sigma _{3}}  - \sigma
_{1} \sigma _{3} \omega _{\sigma _{1} \sigma _{3}}  - \sigma _{2}
\sigma _{3} \omega _{\sigma _{2} \sigma _{3}}  + \sigma _{3}
\sigma _{4} \omega
_{\sigma _{3} \sigma _{4}}}   \\
 { + {\left[ {2\alpha _{3} - (\alpha + \alpha _{3} +1)\sigma _{3}}
\right]}\omega _{\sigma _{3}}  - \alpha _{3} \sigma _{1} \omega
_{\sigma _{1}}  - \alpha _{3} \sigma _{3} \omega _{\sigma _{3}}  +
\alpha _{3} \sigma
_{4} \omega _{\sigma _{4}}  - \alpha \alpha _{3} \omega = 0} \\
 \end{array}} \hfill \\
 {\sigma _{4} \omega _{\sigma _{4} \sigma _{4}}  - \sigma _{1} \omega
_{\sigma _{1} \sigma _{4}}  - \sigma _{2} \omega _{\sigma _{2}
\sigma _{4}} - \sigma _{3} \omega _{\sigma _{3} \sigma _{4}} + (1
- \alpha )\omega _{\sigma _{4}}  + \omega =
0\,\,\,\,\,\,\,\,\,\,\,\,\,\,\,\,\,\,\,\,\,\,\,\,\,\,\,\,\,\,\,\,\,\,\,\,\,\,\,\,\,\,\,\,\,\,\,\,}
\hfill \\
\end{array}}}  \right.}
\end{equation}
Considering the solutions (\ref{e25})-(\ref{e212}) of the system
(\ref{e24}), we define  the solutions $\omega_i(\sigma),
i=1,...,8$  of the system (\ref{eq43}) and substituting those
found solutions into the expression (\ref{eq41}), we get some
fundamental solutions of the equation (\ref{eq11})
\begin{equation}
 \label{eq44} q_1(x,x_0)=k_1\left( r^{2} \right)^{ -
 \alpha}\textrm{H}_{4,3}^0(\alpha,\alpha_1,\alpha_2,\alpha_3;2\alpha_1,2\alpha_2,2\alpha_3;\sigma),
\end{equation}

\begin{equation}
 \label{eq45}\begin{array}{l}q_2(x,x_0)=k_2\left( r^{2} \right)^{ 2\alpha_1- \alpha-1}
\left(x_1x_{01}\right)^{1-2\alpha_1}\\\,\,\,\,\,\,\,\,\,\,\,\,\,\,\,\,\,\,\,\,\,\,\,\,\,\,\,\,\,\,\,\,\,\,\,\,\cdot\textrm{H}_{4,3}^0(1+\alpha-2\alpha_1,1-\alpha_1,\alpha_2,\alpha_3;2-2\alpha_1,2\alpha_2,2\alpha_3;\sigma),\\\end{array}
\end{equation}

\begin{equation}
 \label{eq46} \begin{array}{l} q_3(x,x_0)=k_3\left( r^{2} \right)^{2\alpha_2 - \alpha-1}
\left(x_2x_{02}\right)^{1-2\alpha_2}\\\,\,\,\,\,\,\,\,\,\,\,\,\,\,\,\,\,\,\,\,\,\,\,\,\,\,\,\,\,\,\,\,\,\,\,\,\cdot\textrm{H}_{4,3}^0(1+\alpha-2\alpha_1,\alpha_1,1-\alpha_2,\alpha_3;2\alpha_1,2-2\alpha_2,2\alpha_3;\sigma),\\\end{array}
\end{equation}

\begin{equation}
 \label{eq47}\begin{array}{l}q_4(x,x_0)=k_4\left( r^{2} \right)^{2\alpha_3 - \alpha-1}
\left(x_3x_{03}\right)^{1-2\alpha_3}\\\,\,\,\,\,\,\,\,\,\,\,\,\,\,\,\,\,\,\,\,\,\,\,\,\,\,\,\,\,\,\,\,\,\,\,\,\cdot\textrm{H}_{4,3}^0(1+\alpha-2\alpha_1,\alpha_1,\alpha_2,1-\alpha_3;2\alpha_1,2\alpha_2,2-2\alpha_3;\sigma),\\\end{array}
\end{equation}

 \begin{equation}
 \label{eq48} \begin{array}{l} q_5(x,x_0)=k_5\left( r^{2} \right)^{2\alpha_1+2\alpha_2-\alpha-2}
\left(x_1x_{01}\right)^{1-2\alpha_1}\left(x_2x_{02}\right)^{1-2\alpha_2}\\
\,\,\,\,\,\,\,\,\,\,\,\,\,\,\,\,\,\,\,\,\,\,\,\,\,\,\,\,\,\,\,\,\,\,\,\,\cdot\textrm{H}_{4,3}^0(2+\alpha-2\alpha_1-2\alpha_2,1-\alpha_1,1-\alpha_2,\alpha_3;2-2\alpha_1,2-2\alpha_2,2\alpha_3;\sigma),\\\end{array}
\end{equation}

\begin{equation}
 \label{eq49}  \begin{array}{l} q_6(x,x_0)=k_6\left( r^{2} \right)^{2\alpha_1+2\alpha_3-\alpha-2}
\left(x_1x_{01}\right)^{1-2\alpha_1}\left(x_3x_{03}\right)^{1-2\alpha_3}\\
\,\,\,\,\,\,\,\,\,\,\,\,\,\,\,\,\,\,\,\,\,\,\,\,\,\,\,\,\,\,\,\,\,\,\,\,\cdot\textrm{H}_{4,3}^0(2+\alpha-2\alpha_1-2\alpha_3,1-\alpha_1,\alpha_2,1-\alpha_3;2-2\alpha_1,2\alpha_2,2-2\alpha_3;\sigma),\\\end{array}\end{equation}

\begin{equation}
 \label{eq510} \begin{array}{l}q_7(x,x_0)=k_7\left( r^{2} \right)^{2\alpha_2+2\alpha_3-\alpha-2}
\left(x_2x_{02}\right)^{1-2\alpha_2}\left(x_3x_{03}\right)^{1-2\alpha_3}\\
\,\,\,\,\,\,\,\,\,\,\,\,\,\,\,\,\,\,\,\,\,\,\,\,\,\,\,\,\,\,\,\,\,\,\,\,\cdot\textrm{H}_{4,3}^0(2+\alpha-2\alpha_2-2\alpha_3,\alpha_1,1-\alpha_2,1-\alpha_3;2\alpha_1,2-2\alpha_2,2-2\alpha_3;\sigma),\\\end{array}
\end{equation}

\begin{equation}
 \label{eq511} \begin{array}{l} q_8(x,x_0)=k_8\left( r^{2}
\right)^{2\alpha_1+2\alpha_2+2\alpha_3-\alpha-3}
\left(x_1x_{01}\right)^{1-2\alpha_1}\left(x_2x_{02}\right)^{1-2\alpha_2}\left(x_3x_{03}\right)^{1-2\alpha_3}\\
\,\,\,\,\,\,\,\,\,\,\,\,\,\,\,\,\,\,\cdot\textrm{H}_{4,3}^0(3+\alpha-2\alpha_1-2\alpha_2-2\alpha_3,1-\alpha_1,1-\alpha_2,1-\alpha_3;2-2\alpha_1,2-2\alpha_2,2-2\alpha_3;\sigma),\\\end{array}
\end{equation}
where $k_1,...,k_8$ are constants which will be determined at
solving boundary value problems for equation (\ref{eq11}).

\section{Singularity properties of fundamental solutions}

Let us show that the found solutions (\ref{eq44})-(\ref{eq511})
have a singularity. We choose a solution $q_1(x,x_0)$. For this
aim we use the expansion (\ref{eq18}) for the  confluent
hypergeometric function (\ref{e23}). As a result, solution
(\ref{eq44}) can be written as follows

$$q_{1} \left( {x,x_{0}}  \right) =  r^{2 - p}r_{1}^{ -
2\alpha _{1}} r_{2}^{ - 2\alpha _{2}}  r_{3}^{ - 2\alpha _{3}}
f\left( {r^{2},r_{1}^{2} ,r_{2}^{2} ,r_{3}^{2}}  \right),$$ where
\begin{equation}\label{eq51} \begin{array}{l}
 f\left( {r^{2},r_{1}^{2} ,r_{2}^{2} ,r_{3}^{2}}\right) =k_{1}
{\sum\limits_{l,m,n,s = 0}^{\infty}  {{\sum\limits_{q = 0}^{s}
{{\sum\limits_{I(3,s)} {A\left( {s,q} \right)s!C_{s}^{q} {\frac{{(
- 1)^{q}\left( {\alpha}  \right)_{l + m + n + s} \left( {\alpha
_{1}}  \right)_{i + l + m} \left( {\alpha _{2}}  \right)_{j + l +
n} \left( {\alpha _{3}}  \right)_{k + m + n}}} {{\left( {\alpha}
\right)_{s} \left( {1 - \alpha}  \right)_{q} \left( {2\alpha _{1}}
\right)_{i + l + m} \left( {2\alpha _{2}}  \right)_{j + l + n}
\left( {2\alpha _{3}}  \right)_{k + m + n} i!j!k!n!m!l!q!}}}}}} }}
}\\
 \cdot\left( {1 - {\frac{{r^{2}}}{{r_{1}^{2}}} }} \right)^{i + l + m}\left( {1 -
{\frac{{r^{2}}}{{r_{2}^{2}}} }} \right)^{j + l + n}\left( {1 -
{\frac{{r^{2}}}{{r_{3}^{2}}} }} \right)^{k + m + n}\left(
{{\frac{{1}}{{2}}}\lambda r} \right)^{2q}\\\cdot F\left( {2\alpha
_{1} - \alpha - j - k,\alpha _{1} + i + l + m;2\alpha _{1} + i + l
+ m;1
-{\frac{{r^{2}}}{{r_{1}^{2}}} }} \right) \\
 \cdot F\left( {2\alpha _{2} - \alpha - i - k - m,\alpha _{2} + j + l +
n;2\alpha _{2} + j + l + n;1 - {\frac{{r^{2}}}{{r_{2}^{2}}} }} \right) \\
 \cdot F\left( {2\alpha _{3} - \alpha - i - j - l,\alpha _{3} + k + m +
n;2\alpha _{3} + k + m + n;1 - {\frac{{r^{2}}}{{r_{3}^{2}}} }}
\right){}_{0}F_{1} \left( {1 - \alpha + q;{\frac{{1}}{{4}}}\lambda
^{2}r^{2}} \right), \end{array}\\
\end{equation}
$$k_1=\frac{4^{\alpha_1+\alpha_2+\alpha_3-1}\Gamma(\alpha)\Gamma(\alpha_1)\Gamma(\alpha_2)\Gamma(\alpha_3)}{\pi^{p/2}\Gamma(2\alpha_1)\Gamma(2\alpha_2)\Gamma(2\alpha_3)}.$$
Following the work \cite{HK} and applying several times a
well-known summation formula \cite[p.61,(14)]{EM}
$$F(a,b;c;1)=\frac{\Gamma(c)\Gamma(c-a-b)}{\Gamma(c-a)\Gamma(c-b)},\,\,\,c \neq 0,-1,-2,..., c-a-b>0,$$
it is easy to show that
\begin{equation}\label{eq52}
f\left( {0,r_{10}^{2} ,r_{20}^{2} ,r_{30}^{2}}  \right) =
\frac{4^{\alpha_1+\alpha_2+\alpha_3-1}}{\pi^{p/2}}\Gamma\left(\frac{p-2}{2}\right),
p>2.
\end{equation}

Expressions (\ref{eq51}) and (\ref{eq52}) give us the possibility
to conclude that the solution $q_1(x,x_0)$ reduces to infinity of
the order $r^{2-p}$ at $r \to 0$. Similarly it is possible to be
convinced that solutions $q_i(x,x_0),\,\,i=2,3,...,8$ also reduce
to infinity of the order $r^{2-p}$  when $r \to 0$.

\end{document}